\documentclass[12pt]{article}

\usepackage{fontenc}
\usepackage{amsmath}
\usepackage{amsfonts}
\usepackage{delarray}
\usepackage{theorem}

\newtheorem{Thm}{Theorem}

\newtheorem{Lem}{Lemma}

\def\C{\mathbb C}
\def\Z{\mathbb Z}
\def\N{\mathbb N}

\begin{document}

\title{On the algebraic solutions of the sixth Painlev\'e   equation related to second order Picard-Fuchs
equations}
\author{
{Bassem \textsc{Ben Hamed\dag}}\hspace{0.5cm}{Lubomir \textsc{Gavrilov\ddag}}}
\maketitle

\hspace*{0.5cm}\dag \hspace{0.2cm}Facult\'e des Sciences de Sfax, D\'epartement de
Math\'ematiques\\
\hspace*{2.5cm}BP 802, route Soukra km 4 Sfax 3018, Tunisie\\
\hspace*{0.5cm}\ddag \hspace{0.2cm}Laboratoire Emile Picard, CNRS UMR 5580, Universit\'e Paul
Sabatier\\
\hspace*{2cm}118, route de Narbonne, 31062 Toulouse Cedex 4, France\\

\begin{abstract}
We describe  two algebraic solutions of the sixth Painlev\'e
equation which are related to (isomonodromic) deformations of
Picard-Fuchs equations of order two.
\end{abstract}

\section{Statement of the result}
In this note we describe  two algebraic solutions of the following
Painlev\'e VI ( $\mathcal P_{VI}$) equation
\begin{eqnarray*}
 \frac{d^2 \lambda}{dt^2} &= &\frac{1}{2}(\frac{1}{\lambda}+ \frac{1}
{\lambda - 1} + \frac{1}{\lambda - t})(\frac{d\lambda}{dt})^2 - (\frac{1}{t} + \frac{1}{t - 1}
 + \frac{1}{\lambda - t})\frac{d\lambda}{dt}\\
& &+ \frac{\lambda (\lambda - 1)(\lambda - t)}
 {t^2 (t^2 - 1)} [ \alpha + \beta \frac{t}{\lambda ^2} + \gamma \frac{t - 1}{(\lambda - 1)^2}
 +  \delta \frac{t(t - 1)}{(\lambda - t)^2} ]  .
\end{eqnarray*}
related to deformations of Picard-Fuchs equations of special type.
Recall that the $\mathcal P_{VI}$ equation governs the
isomonodromic deformations of the second order Fuchsian equations
\begin{equation}
x''+p_1(s)x'+p_2(s)x=0,\quad  '= \frac{d}{ds} \label{fuchs}
\end{equation}
with 5 singular points, one of which is apparent \cite{8}. Suppose
that the solution  of the   Fuchsian equation (\ref{fuchs}) is
given by an Abelian integral
$$
x(s)= \int_{\gamma_s} \omega
$$
where $\omega$ is a rational one-form on $\mathbb{C}^2$, $\Gamma_s
\subset \mathbb{C}^2$ is a family of algebraic curves depending
rationally on $s$, and $\gamma_s \subset \Gamma_s$ is a continuous
family of closed loops. Then the equation (\ref{fuchs}) is said to
be of Picard-Fuchs type and its monodromy group is conjugated to a
subgroup of $\mathbf{Gl}_2(\overline{\mathbb{Q}})$ (generically
$\mathbf{Gl}_2(\mathbb{Z})$). For this reason any continuous
deformation
$$
a \rightarrow \Gamma_{s,a}
$$
of the family $\Gamma_s$ induces an isomonodromic deformation of
(\ref{fuchs}). If in addition $\Gamma_{s,a}$ depends algebraically
in $a$, the coefficients of (\ref{fuchs}) are also algebraic
functions in $a$, and hence they provide an algebraic solution of
$\mathcal P_{VI}$.

From now on we denote
$$
(\alpha_0,\alpha_1,\alpha_2,\alpha_3)=
 (\alpha,-\beta,\gamma,\frac{1}{2} - \delta).
 $$

%
%
%
%
%
Our main result is the following
\begin{Thm}
\label{t7}
The pencil of  $\mathcal P_{VI}(\alpha )$ equations

\begin{equation}\label{pencil}
(\alpha_0,\alpha_1,\alpha_2,\alpha_3) =
(\frac{1}{8},\frac{s}{8},\frac{s}{8},\frac{s}{8}), s\in \C
\end{equation}
has a common algebraic solution parameterized as
\begin{equation}
\lambda =\frac{a^2(2-a)}{a^2-a+1}, t=\frac{a^3(2-a)}{2a-1}, a\in \C .
\label{eq1}
\end{equation}
The $\mathcal P_{VI}(\alpha )$ equation with
$$
\alpha = (\frac{1}{8}, \frac{1}{2},0,0)
$$
has an algebraic solution parameterized as
\begin{equation}
\lambda =\frac {a(a-2)(2a^2+a+2)}{a^2-7a+1}, t=\frac{a^3(2-a)}{2a-1}, a\in \C
\label{eq2}
\end{equation}
\end{Thm}
The meaning of these solutions is the following. Consider the family of elliptic curves
\begin{equation}
\Gamma _s = \{(\xi ,\eta )\in \C^2 :
\eta ^2 + \frac{3}{2a-1} \xi ^4 - \frac{4(a+1)}{2a-1} \xi ^3 + \frac{6a}{2a-1} \xi ^2= s \}
\label{curve}
\end{equation}
and let $\gamma (s) \in H_1(\Gamma _s,\Z)$ be a family of cycles
depending continuously on $s\in \C$. The Abelian integral of first
kind
$$
\int_{\gamma (s)} \frac{d \xi}{\eta }
$$
satisfies a Picard-Fuchs equation of second order   depending on a
parameter $a$, defining an isomonodromy deformation of the
equation. This deformation corresponds then to an algebraic
solution of $\mathcal P_{VI}(\alpha )$  given by (\ref{eq1}). In a
similar way, the Abelian integral of second kind
$$
\int_{\gamma (s)} \frac{(3\xi^2- 2(a+1))\xi d\xi }{\eta }
$$
satisfies a Picard-Fuchs equation of second order. The
isomonodromy deformation of this equation with respect    to $a$
is described by the solution (\ref{eq2}) of $\mathcal P_{VI}$
equation.

Algebraic solutions of $\mathcal P_{VI}$ were found by many
authors, e.g. Hitchin\cite{7}, Manin\cite{10},
Dubrovin-Mazzocco\cite{dub00}, Boalch\cite{boalch}. Dubrovin and
Mazzocco classified all algebraic solutions of the $\mathcal
P_{VI}$ equation corresponding to
$$
(\alpha_0,\alpha_1,\alpha_2,\alpha_3) = ( \frac{1}{2} (2 \mu
-1)^2,0,0,0), \mu \in \mathbb{R} .
$$
It turns out that these solutions, up to symmetries, are in a
one-to-one correspondence with the regular polyhedra in the three
dimensional space.
Our solution (\ref{eq1}) with
$$(\alpha_0,\alpha_1,\alpha_2,\alpha_3) =(1/8,0,0,0)$$
 corresponds then to the tetrahedron solution of Dubrovin-Mazzocco
 ($\mu =+1/4$ ). It is identified to their solution
 $(A_3)$ via the Okamoto type transformation
 (1.24),(1.25), see \cite{dub00}. It is remarkable that the same  solution, but for
$$(\alpha_0,\alpha_1,\alpha_2,\alpha_3) =(1/8,1/8,1/8,1/8)$$
was also found by Hitchin\cite{7}. This shows that (\ref{eq1}) is
a common solution to the family (\ref{pencil}) of $\mathcal
P_{VI}$ equations. It is clear that for transcendental values of
$s$ in (\ref{pencil}) the corresponding isomonodromic family of
Fuchs equations (\ref{fuchs}) can not be of Picard-Fuchs type.

 The paper is
organized as follows. In the next section we recall briefly,
following \cite{8}, the relationship between $\mathcal P_{VI}$ and
the isomonodromic deformations of Fuchs equations. In section
\ref{pf} we deduce the relevant Picard-Fuchs equation and
establish the main result.

The present text is an abridged version of \cite{bassem}.

\section{The Garnier system and the $\mathcal P_{VI}$ equation.}
Consider a Fuchsian differential equation
$$
x''+p_1(s)x'+p_2(s)x=0,\quad  '= \frac{d}{ds} \eqno{(\ref{fuchs})}
$$
with {\it five singular points, exactly one  of which  is
apparent}. After a bi-rational change of the independent variable
$s$ and a linear change of the dependent variable $x$ (involving
$s$) we may suppose that the singular points are $0,1,t,\lambda
,\infty$, where the singularity $\lambda $ is apparent and the
corresponding Riemann scheme is
$$
\begin{array}({ccccc})
0 &1 & t &\lambda & \infty \\
0 &0 & 0&               0& \alpha \\
\theta_1 & \theta _2&\theta _3& k & \alpha +\theta_{\infty}
\end{array}, n\in \N, 2\alpha + \sum_i \theta _i +n = 3 .
$$
In what follows we shall always suppose that $n=2$  (which is
satisfied generically)).

The coefficients $p_1,p_2$ are easily computed to be
$$
p_1(s)=\frac{1-\theta _1}{s}+\frac{1-\theta _2}{s-1}+\frac{1-\theta _3}{s-t}
-\frac{1}{t-\lambda}
$$
$$
p_2(s)=
\frac{k}{s(s-1)}- \frac{t(t-1)K}{s(s-1)(s-t)}+\frac{\lambda
(\lambda-1)\mu}{s(s-1)(s-\lambda)}
$$
where $\mu$ is a  constant
$$
k=\frac{1}{4}\{(\sum_{i=1}^{3} \theta_i - 1)^2-\theta_{\infty}^2\} .
$$
The compatibility condition for the
singular point $\lambda $ to be apparent reads
\begin{eqnarray*}
K=K(\lambda ,\mu ,t) &=& \frac{1}{t(t-1)}[\lambda (\lambda -1)(\lambda -t)\mu^2 -
\{\theta_2 (
\lambda -1)(\lambda -t)\\
& & +\theta_3\lambda(\lambda -t) + (\theta_1 -1)\lambda (\lambda
-1)\}\mu +k\lambda]  .
\end{eqnarray*}
 From the discussion above it is seen that the Fuchs equation (\ref{fuchs}) depends on the
parameters $\theta_0, \theta _1, \theta _t, \theta _\infty, \lambda , \mu ,t$. Let us denote this
equation by $E_\theta (\lambda ,\mu ,t)$.

\begin{Thm}
\label{garnier}
$\lambda (t),\mu (t)$ is a solution of the Garnier system
\begin{eqnarray*}
\frac{d\lambda}{dt} & = &\frac{\partial K}{\partial \mu} \\
\frac{d\mu}{dt} & = &-\frac{\partial K}{\partial \lambda}  .
\end{eqnarray*}
if and only if
the induced deformation of $E_\theta (\lambda ,\mu ,t)$ is isomonodromic.
\end{Thm}
It is straightforward to check that the sixth Painlev\'e system
$\mathcal P_{VI}(\alpha )$ with parameters
\begin{equation}
\alpha =(\frac{1}{2}\theta _{\infty}^2 , \frac{1}{2}\theta _0^2,\frac{1}{2}\theta
_1^2,\frac{1}{2}\theta _t^2)
\label{thetas}
\end{equation}
is equivalent to the  Garnier system. We get therefore the following\\
{\bf Corollary.} If
$$(t,\lambda ,\mu )\rightarrow (t,\lambda (t),\mu (t))
$$
is an isomonodromic deformation of $E_\theta (\lambda ,\mu ,t)$, then $\lambda (t)$ is a
solution of $\mathcal P_{VI}(\alpha )$ equations with parameters given by (\ref{thetas}).

\section{Picard-Fuchs equations}
\label{pf}

 In this section we restrict our attention to the
deformation
$$
f_a(\xi ,\eta) = \eta ^2 + \frac{3}{2a-1} \xi ^4 - \frac{4(a+1)}{2a-1} \xi ^3 + \frac{6a}{2a-1} \xi
^2, a\in \C
$$
 of the singularity
$\eta ^2 +\xi^4$
 of
type $A_3$, see \cite{3}. The critical values of $f_a(\xi ,\eta)$
are
$$0,1, t=\frac{a^3(2-a)}{2a-1} .$$
Consider the locally trivial smooth fibration
$$   f^{-1}(\C \setminus \{0,1,t\}) \rightarrow \C \setminus \{0,1,t\}                    $$
whose fibers the
affine curves $\Gamma _s$, (\ref{curve}), $s\in \C \setminus \{0,1,t\}$. Each $\Gamma _s$ is
topologically a torus with two removed points. Hence
$\dim H_1(\Gamma _s,\Z)= \dim H^1_{DR}(\Gamma _s,\C)=3$. Therefore if $\gamma (s) \in
H_1(\Gamma _s,\Z)$ is a family of cycles depending continuously on $s$, then the Abelian
integral
$$
I(s)= \int_{\gamma (s)} \omega , \omega = P(\xi ,\eta) d\xi  + Q(\xi ,\eta) d\eta, P,Q \in \C[\xi ,
\eta]
$$
satisfies a Fuchsian differential equation of order three, whose
coefficients are polynomials in $s,a$. In the case when the
differential form $\omega $ has no residues, it satisfies a second
order equation. Explicitly, if $\gamma _1(s)$, $\gamma _2(s)$, is
a continuous family of cycles generating the homology group of the
compactified elliptic curve $\Gamma _s$, then the equation reads
$$
\det \left(
\begin{array}{ccc}
x& x'& x'' \\
\int_{\gamma _1(s)} \omega & (\int_{\gamma _1(s)}  \omega )'&(\int_{\gamma _1(s)} \omega
)''\\
\int_{\gamma _2(s)} \omega & (\int_{\gamma _2(s)}  \omega )'&(\int_{\gamma _2(s)} \omega
)''
\end{array}\right)= 0 .
$$
It follows from the Picard-Lefschetz formula and the moderate growth of the integrals, that
the coefficients of the above differential equations are rational in
 $s,a$.     A local analysis of the singularities shows for instance that
$$
\det \left(
\begin{array}{cc}
\int_{\gamma _1(s)} \omega & (\int_{\gamma _1(s)}  \omega )'\\
\int_{\gamma _2(s)} \omega & (\int_{\gamma _2(s)}  \omega )'
\end{array} \right) =
\frac{p(s,a)}{s(s-1)(s-t)}
$$
where $p(s,a)$ is a polynomial in $s,a$. If we put $\omega =dx/y$ then $\int_{\gamma _1(s)}
\omega$ grows no faster than $s^{1/4-1/2}$ at $\infty$ (for a fixed $a$). Thus
$$
\frac{p(s,a)}{s(s-1)(s-t)}
$$
grows at infinity no faster than $s^{-1/2-1}$ and hence no faster
than $s^{-2}$. It is expected therefore than $p(s,a)$ is of degree
one in $s$ and the corresponding root, which we denote by $\lambda
$, is an apparent singularity for the Picard-Fuchs equation in
consideration. We are therefore in a position to apply Theorem
\ref{garnier},  provided that the deformation of the Fuchs
equation with respect to the parameter $a$ is isomonodromical.
Indeed, the monodromy group of our equation is contained in
$SL(2,\Z)$ which shows that any deformation of this equation is
isomonodromical. The Picard-Lefschetz formula shows that the
monodromy group in question is generated, up to conjugacy,  by the
matrices
\begin{equation}
\begin{array}({cc})
1 & 1\\
0& 1
\end{array},\quad
\begin{array}({cc})
1 & 0\\
1& 1
\end{array} .
\label{monodromy}
\end{equation}

 To deduce an explicit formula for the corresponding algebraic solution of  $\mathcal P_{VI}$ we need
 explicit formulae for the Picard-Fuchs equations.
\begin{Lem}
Let $\gamma (s)\in H_1(\Gamma _s,\Z)$ be a family of cycles
depending continuously on $s$. The complete elliptic integrals of
first and second kind
$$
x(s)= \int_{\gamma (s)} \frac{d\xi }{\eta },\;\;\;
y(s)=\int_{\gamma (s)} \frac{(3\xi^2- 2(a+1))\xi d\xi }{\eta }$$
satisfy Picard-Fuchs equations of the form
$$
a_0(s) x'' + a_1(s) x' + a_2(s) x= 0$$
 $$
 b_0(s) y'' + b_1(s) y' +
b_2(s) y= 0
$$
where
\begin{eqnarray*}
a_0(s)&=&s(s-1)((2a-1)s+a^3(a-2))((a^2-a+1)s+a^2(a-2)) \\
a_1(s)&=&2(2a-1)(a^2-a+1)s^3+(a^6-3a^5+9a^4-19a^3+9a^2-3a+1)s^2 \\
                 & &+2a^2(a-2)(a^4-2a^3-2a+1)s-a^5(a-2)^2 \\
a_2(s)&=&(2a-1)[27(a^2-a+1)s^2-(a-2)(2a^4-a^3-60a^2-a+2)s \\
& & +a^2(a-2)^2(10a^2+11a+10)]/144 \\
b_0(s)&=&s(s-1)((2a-1)s+a^3(a-2))((a^2-7a+1)s-a(a-2)(2a^2+a+2)) \\
b_1(s)&=&(2a-1)s[(a^2-7a+1)s^2-2a(a-2)(2a^2+a+2)s \\
& & -a(a-2)^2(a^4+a^3+a^2+a+1)] \\
b_2(s)&=& -(2a-1)[9(a^2-7a+1)s^2-(a-2)(10a^4+31a^3-12a^2+31a+10)s \\
 & & -a(a-2)^2(2a^2+a+2)^2]/144
\end{eqnarray*}
\end{Lem}
The proof of the above Lemma is straightforward, see for instance \cite{dum01}.
It is seen that the roots of $a_0(s)$ are $0,1$ and
$$
\lambda =\frac{a^2(2-a)}{a^2-a+1}, t=\frac{a^3(2-a)}{2a-1}
$$
which implies the algebraic solution (\ref{eq1}). In the same way the roots of $b_0(s)$ provide
the solution (\ref{eq2}).
 The Riemann schemes of the Picard-Fuchs equations for $x(s),y(s)$ are given by
$$
\begin{array}({ccccc})
0 &1 & \frac{a^3(2-a)}{2a-1} &\frac{a^2(2-a)}{a^2-a+1} & \infty \\
0 &0&0&0&  \frac{1}{4} \\
0 &0&0&2&  \frac{3}{4}
\end{array}
$$
and
$$
\begin{array}({ccccc})
0 &1 & \frac{a^3(2-a)}{2a-1} &\frac{a(a-2)(2a^2+a+2)}{a^2-7a+1} & \infty \\
0 &0&0&0&  \frac{1}{4} \\
1 &0&0&2&  -\frac{1}{4}
\end{array}.
$$
The Corollary after Theorem \ref{garnier} implies that the  curve
(\ref{eq1}) is an integral curve of the $\mathcal P_{VI}(\alpha )$
equation with parameters  $ (\alpha_0,\alpha_1,\alpha_2,\alpha_3)
= (\frac{1}{8},0,0,0) $, see (\ref{thetas}). Similarly, the
Fuchsian equation satisfied by the complete elliptic integral of
second kind $y(s)$ provides the algebraic solution (\ref{eq2})  of
$\mathcal P_{VI}(\alpha )$ with $\alpha =
(\frac{1}{8},\frac{1}{2},0,0)$.

It is remarkable that  (\ref{eq1}) was found to be a solution of
$\mathcal P_{VI}(\alpha )$ by Hitchin \cite{7}[p.177], but for
$\alpha = (\frac{1}{8},\frac{1}{8},\frac{1}{8},\frac{1}{8})$ .
After taking the difference  between these two equations we obtain
the following affine equation of the integral curve (\ref{eq1})
$$
-  \frac{t}{\lambda ^2} +  \frac{t - 1}{(\lambda - 1)^2}
 - \frac{t(t - 1)}{(\lambda - t)^2} =0 .
$$
This also shows that (\ref{eq1}) is a common algebraic solution of
the pencil of  $\mathcal P_{VI}(\alpha )$ equations
$$
\alpha  = (\frac{1}{8},\frac{s}{8},\frac{s}{8},\frac{s}{8}), s\in \C .
$$
This completes the proof of Theorem \ref{t7}.

\end{document}